\documentclass[12pt,leqno]{amsart}
\usepackage{amsmath,amssymb,amsfonts}
\usepackage{hyperref}

\newtheorem{theorem}{Theorem}

\theoremstyle{definition}
\newtheorem{problem}[theorem]{Problem}
\newtheorem{example}[theorem]{Example}

\begin{document}

\title{Two Hamiltonian Cycles}

\author{Vaidy Sivaraman}\address{Department of Mathematics and Statistics, Mississippi State University, Mississippi State, MS 39762, USA}\email{vs692@msstate.edu}
\author{Thomas Zaslavsky}\address{Dept.\ of Mathematical Sciences, Binghamton University (SUNY), Binghamton, NY 13902-6000, USA}\email{zaslav@math.binghamton.edu}

\begin{abstract}
If the line graph of a graph $G$ decomposes into Hamiltonian cycles, what is $G$?  We answer this question for decomposition into two cycles.
\end{abstract}



\date{\today}

\maketitle

\section*{Prelude}

Two problems that are known to be NP-complete are:  Is a cubic graph Hamiltonian \cite{gjt}?  Does a graph $L$ decompose into two Hamiltonian cycles \cite{peroche}?  
An example for $L$ is the line graph of a cubic graph $G$: $L(G)$ has a Hamilton decomposition if and only if $G$, the root graph, is Hamiltonian; this was discovered by Kotzig \cite{kotzig} and later Martin \cite{martin}.  
(We say a graph is \emph{Hamilton decomposable} if its edge set partitions into Hamiltonian cycles.)  
The proofs are by polynomial-time constructions.  Since Hamiltonicity of a cubic graph is NP-complete, it follows that Hamiltonian decomposability of a 4-regular graph remains NP-complete when restricted to line graphs of cubic graphs.

An opposite question seems not to have been asked.  (It was briefly considered by the second author for an undergraduate graph theory examination.)

\begin{problem}\label{2hc}
The line graph $L(G)$ of a simple graph $G$ is decomposable into two Hamiltonian cycles.  What is $G$?
\end{problem}

We answer this question by characterizing all root graphs $G$. 
One would like to describe all Hamilton decomposable line graphs.  We offer some thoughts at the end of this paper.

\section*{Presto}

\begin{theorem}\label{t:2hc}
The graph $G$ of Problem \ref{2hc} is either $K_{1,5}$, or the first subdivision of a $4$-regular graph $G'$ that decomposes into two Hamiltonian cycles, or a Hamiltonian cubic graph.
\end{theorem}

The first subdivision $SG$ of a graph $G$ is obtained by subdividing every edge into a path of length 2.  The 4-regular graph $G'$ is not assumed to be simple.

The first case gives line graph $K_5$, which is obviously a solution to Problem \ref{2hc}.  In the second case $G=SG'$ is not itself Hamiltonian; we shall explain how its line graph decomposes in terms of a known Hamilton decomposition of $G'$.  If $G'$ has $n$ vertices, there are $2^n$ different ways to decompose $L(G)$ using that decomposition of $G'$.  In the third case we assume a Hamiltonian cycle is known.  If there are $2n$ vertices in $G$, there are $2^n$ different ways to decompose $L(G)$ using that Hamiltonian cycle.

The reductions from the line graph and its Hamiltonian cycles to the root graph and its Hamiltonian cycle(s) can be performed in polynomial time, since a graph is reconstructible from its line graph in polynomial time \cite{lehot}.

\section*{Moderato}

We describe the construction of the Hamilton decomposition of our 4-regular $L(G)$ in each case other than the familiar $K_5$.  We give the cubic construction, from \cite{kotzig, martin}, for completeness.

\begin{example}\label{x:4reg}
This example is more general than is needed for Theorem \ref{t:2hc}.  
Suppose $G$ is the first subdivision of a $2h$-regular graph $G'$, $h\geq2$, that decomposes into $h$ Hamiltonian cycles $H_1',\ldots,H_h'$.  We convert $L(SH_i')$ into a Hamiltonian cycle $H_i$ of $L(G)$ as follows.

A divalent vertex in $G$, incident with edges, say, $c$ and $d$ derived by subdividing $c'$ in $G'$, gives an edge $cd$ that we put in $H_i$ if $c' \in H_i'$.  

Consider a $2h$-valent vertex $u$ in $G$.  At $u$ there are two edges of $G$ belonging to each $SH_i'$, call them $e_i$ and $f_i$.  The edges $e_1,f_1,\ldots,e_h,f_h$ become vertices of a $2h$-clique $K_{2h}(u)$ in $L(G)$.  
The graph $K_{2h}(u)$ decomposes into $h$ Hamiltonian paths $H_i(u)$ with distinct endpoints that can be chosen arbitrarily \cite[proof of Theorem 2.3.3]{pearls}; we choose them to have respective endpoints $e_i,f_i$ for each $i=1,\ldots,h$.
Replace the edge $e_if_i$ in $L(SH_i')$ by the path $H_i(u)$.  Doing this for all $2h$-valent vertices in $G$ gives $h$ edge-disjoint 2-factors of $L(G)$:  $H_i$ derived from $L(SH_i')$ for each $i$.  Each $H_i$ is 2-regular by construction and it is connected because $H_i'$ is connected; therefore it is a Hamiltonian cycle of $L(G)$.  We have a decomposition because each edge of $L(G)$ has been assigned to exactly one $H_i$.
\end{example}

In \cite{hdecomp} Muthusamy and Paulraja proved half the cases of the similar result that $L(G')$ itself decomposes into Hamiltonian cycles, and indeed more generally that if $G'$ is any Hamilton decomposable graph then $L(G')$ is also Hamilton decomposable.  Example \ref{x:4reg} differs in that we subdivide $G'$ first, which makes the root graph no longer Hamilton decomposable.

\begin{example}\label{x:cubic}
Suppose $G$ is a cubic graph with a Hamiltonian cycle $H$.  Let $M$ be the complementary matching:  $M = E(G) \setminus E(H)$.  Orient the edges of $M$ arbitrarily.  Expand $H_0 = L(H)$ into two Hamiltonian cycles of $L(G)$, as follows.  

Begin with two copies of $H_0$, labeled $H_1$ and $H_2$.  Consider an edge $ef$ of $L(H)$ corresponding to two edges of $H$ with common vertex $u$ and let $m=uu'$ be the matching edge at $u$.  The vertex $u'$ belongs to two edges of $H$, say $e'$ and $f'$, forming the edge $e'f'$ in $H_0$.  
If $u$ is the tail of $uu'$, replace the edge $ef$ of $H_1$ by the path $emf$ and replace the edge $e'f'$ in $H_2$  by $e'mf'$, but in $H_2$ retain the edge $ef$ and in $H_1$ retain the edge $e'f'$.  Doing this for all matching edges makes $H_1$ and $H_2$ into edge-disjoint Hamiltonian cycles of $L(G)$.
\end{example}

The decompositions depend on the choice of Hamilton decomposition in Example \ref{x:4reg} and of Hamiltonian cycle in Example \ref{x:cubic}.  Thus, in general there may be a great many Hamilton decompositions of $L(G)$ besides the $2^n$ mentioned earlier.

\section*{Andante}

We prove that the examples in Theorem \ref{t:2hc} are the only ones.  We assume given a line graph $L=L(G)$, but not the root graph $G$, and a decomposition of $L$ into Hamiltonian cycles $H_1$ and $H_2$.  A root graph can be calculated quickly from $L$ \cite{lehot}, thus the vertex cliques in $L$ are known.  In particular, we know when $G=K_{1,5}$ because then $L$ is a 5-clique.  If it is not, then we know when $G$ is $(4,2)$-biregular by the existence of 4-cliques in $L$, and $G$ is cubic in the remaining case.  Thus, we can consider each type of root graph separately.  We leave the case $L=K_5$ to the reader.

\emph{Suppose $G$ is $(4,2)$-biregular.}  We can reconstruct $G''$, a graph isomorphic to $G'$, directly from the 4-cliques in $L$.  The vertices $k_i$ of $G''$ are the 4-cliques in $L$, and there is an edge $k_ik_j$ of $G''$ for each edge between a vertex of $L$ in $k_i $ and one in $k_j$.  We know $G'' \cong G'$ because each vertex of $G''$ is the vertex clique of $L=L(G)$ that corresponds to a quadrivalent vertex in $SG'$, that is, a vertex of $G'$.

The graph $G''$ enables us to deduce the two Hamiltonian cycles $H_i''$ of $G''$ from the $H_i$ in $L$.  Since a vertex in $G''$, considered as a 4-clique in $L$, is quadrivalent, each Hamiltonian cycle $H_i$ of $L$ must enter and leave that 4-clique exactly once using separate edges from those used by the other one.  Thus, $H_i$ acts as a divalent subgraph $H_i''$ of $G''$.  This $H_i''$ is connected because $H_i$ is, hence it is a Hamiltonian cycle of $G''$.  That proves the characterization of $G$ in the $(4,2)$-biregular case of $L$.

Note that, if we begin with $G'$ and its Hamiltonian cycles $H_i'$, construct $L = L(SG')$ and its Hamiltonian cycles $H_i'$ as in Example \ref{x:4reg}, and then construct $H_i''$ in $G''$, then the natural isomorphism of $G''$ to $G$ carries each $H_i''$ to $H_i$.

\emph{Suppose $G$ is cubic.}  We reconstruct the Hamiltonian cycle in $G$ from $L$ by examining the vertex triangles in $L$.  No such triangle can have all three edges in a single one of $H_1$ and $H_2$, so there are two kinds of vertex triangle:  type 1 has two edges from $H_1$ and type 2 has two edges from $H_2$.  Suppose $e_1,e=e_2,e_3$ are the vertices of a type 1 vertex triangle in $L$, derived from a vertex $v$ of $G$, let $u_1,w,u_3$ be the other endpoints of $e_1,e,e_3$ in $G$, and say the path $e_1ee_3$ is in $H_1$, so $H_2$ contains the path $e_1e_3$.  The trace of $H_1$ in $G$ follows the re-entrant trajectory $u_1vwvu_3$ (in some direction) and that of $H_2$ follows the simple trajectory $u_1vu_2$.  
If we remove the re-entrance from the trace of $H_1$, then $H_1$ and $H_2$ are identical.

The other endpoint of $e$ in $G$, namely $w$, forms its own vertex triangle in $L$ with vertices $f_1,e,f_3$.  Let $x_1,x_3$ be the other endpoints of $f_1,f_3$ in $G$.  In $L$, $e$ is incident with the four edges $ee_1,ee_3,ef_1,ef_3$.  The first two belong to $H_1$, so the second pair must belong to $H_2$.  Therefore, the trace of $H_2$ in $G$ contains the re-entrant trajectory $x_1wuwx_3$ and that of $H_1$ contains the path $x_1wx_3$.  The re-entrant edge here is the same edge $e=uw$ as in the trace $u_1vwvu_3$ of $H_1$ from the $v$-vertex triangle.  This implies that the re-entrant edges are a well-defined set and a 1-factor of $G$.  The complementary 2-factor of $G$ is the simplified trace of both $H_1$ and $H_2$ (that is, the trace after removing re-entrant edges), thus it is connected, so it is a Hamiltonian cycle of $G$.  This proves the characterization of $G$ in the 4-regular case of $L$.

It is again easy to see that if we begin with cubic $G$ having a Hamiltonian cycle $H$, apply the construction of $L(G)$ and its Hamilton decomposition $H_1$ and $H_2$ into as in Example \ref{x:cubic}, and return to $G$ by the preceding construction of a new Hamiltonian cycle derived from $H_1$ and $H_2$, the new Hamiltonian cycle of $G$ is nothing but $H$ again.

\section*{Coda}

One naturally thinks about characterizing all line graphs that are Hamilton decomposable, say into $h$ Hamiltonian cycles.  The root graph is then $(k+1,2h-k+1)$-biregular for some $k\leq h$.  The star, the case $k=0$, generalizes to $K_{2h+1}$, whose line graph $K_{2h+1}$ decomposes into $h$ Hamiltonian cycles.  The 4-regular solution generalizes to the case $k=1$, where we have a $(2h,2)$-biregular root graph as in Example \ref{x:4reg}.  In the case $k=h$ we have the line graph of an $h+1$-regular graph $G$, generalizing our cubic root graph.  If $G$ can be proved Hamilton decomposable (for odd $h$) or Hamilton decomposable aside from a 1-factor (for even $h$), then this case is solved by \cite{hdecomp}.  If $G$ can be proved merely to be Hamiltonian or to have a Hamiltonian 3-factor, this same case is solved by \cite{bhms}.  Aside from that, there is the virtual certainty that each value of $k$ introduces a new family of solutions.


\end{document}